\documentclass[12pt]{amsart}

\usepackage{amsmath,amsfonts,amssymb,amscd,verbatim,delarray,verbatim}
\newcommand{\F}{{\mathbb F}}

\newcommand{\Q}{{\mathbb Q}}

\begin{document}
\numberwithin{equation}{section}

\newtheorem{theorem}{Theorem}[section]
\newtheorem{lemma}[theorem]{Lemma}
\newtheorem{prop}[theorem]{Proposition}
\newtheorem{proposition}[theorem]{Proposition}
\newtheorem{corollary}[theorem]{Corollary}
\newtheorem{corol}[theorem]{Corollary}
\newtheorem{conj}[theorem]{Conjecture}

\theoremstyle{definition}
\newtheorem{defn}[theorem]{Definition}
\newtheorem{example}[theorem]{Example}
\newtheorem{examples}[theorem]{Examples}
\newtheorem{remarks}[theorem]{Remarks}
\newtheorem{remark}[theorem]{Remark}
\newtheorem{algorithm}[theorem]{Algorithm}
\newtheorem{quest}[theorem]{Question}
\newtheorem{problem}[theorem]{Problem}
\newtheorem{subsec}[theorem]{}


\def\toeq{{\stackrel{\sim}{\longrightarrow}}}
\def\into{{\hookrightarrow}}


\def\alp{{\alpha}}  \def\bet{{\beta}} \def\gam{{\gamma}}
 \def\del{{\delta}}
\def\eps{{\varepsilon}}
\def\kap{{\kappa}}                   \def\Chi{\text{X}}
\def\lam{{\lambda}}
 \def\sig{{\sigma}}  \def\vphi{{\varphi}} \def\om{{\omega}}
\def\Gam{{\Gamma}}   \def\Del{{\Delta}}
\def\Sig{{\Sigma}}   \def\Om{{\Omega}}
\def\ups{{\upsilon}}


\def\F{{\mathbb{F}}}
\def\BF{{\mathbb{F}}}
\def\BN{{\mathbb{N}}}
\def\Q{{\mathbb{Q}}}
\def\Ql{{\overline{\Q }_{\ell }}}
\def\CC{{\mathbb{C}}}
\def\R{{\mathbb R}}
\def\V{{\mathbf V}}
\def\D{{\mathbf D}}
\def\BZ{{\mathbb Z}}

\def\XX{\mathbf{X}^*}
\def\xx{\mathbf{X}_*}

\def\AA{\Bbb A}
\def\BA{\mathbb A}
\def\HH{\mathbb H}
\def\PP{\Bbb P}

\def\Gm{{{\mathbb G}_{\textrm{m}}}}
\def\Gmk{{{\mathbb G}_{\textrm m,k}}}
\def\GmL{{\mathbb G_{{\textrm m},L}}}
\def\Ga{{{\mathbb G}_a}}

\def\Fb{{\overline{\F }}}
\def\Kb{{\overline K}}
\def\Yb{{\overline Y}}
\def\Xb{{\overline X}}
\def\Tb{{\overline T}}
\def\Bb{{\overline B}}
\def\Gb{{\bar{G}}}
\def\Ub{{\overline U}}
\def\Vb{{\overline V}}
\def\Hb{{\bar{H}}}
\def\kb{{\bar{k}}}

\def\Th{{\hat T}}
\def\Bh{{\hat B}}
\def\Gh{{\hat G}}

\def\cF{{\mathfrak{F}}}
\def\cC{{\mathcal C}}
\def\cU{{\mathcal U}}

\def\Xt{{\widetilde X}}
\def\Gt{{\widetilde G}}

\def\gg{{\mathfrak g}}
\def\hh{{\mathfrak h}}
\def\lie{\mathfrak a}

\def\textrm#1{\text{\textnormal{#1}}}

\def\GL{\textrm{GL}}            \def\Stab{\textrm{Stab}}
\def\Gal{\textrm{Gal}}          \def\Aut{\textrm{Aut\,}}
\def\Lie{\textrm{Lie\,}}        \def\Ext{\textrm{Ext}}
\def\PSL{\textrm{PSL}}          \def\SL{\textrm{SL}}
\def\loc{\textrm{loc}}
\def\coker{\textrm{coker\,}}    \def\Hom{\textrm{Hom}}
\def\im{\textrm{im\,}}           \def\int{\textrm{int}}
\def\inv{\textrm{inv}}           \def\can{\textrm{can}}
\def\id{\textrm{id}}              \def\char{\textrm{char}}
\def\Cl{\textrm{Cl}}
\def\Sz{\textrm{Sz}}
\def\ad{\textrm{ad\,}}
\def\SU{\textrm{SU}}
\def\PSL{\textrm{PSL}}
\def\PSU{\textrm{PSU}}
\def\rk{\textrm{rk}}
\def\PGL{\textrm{PGL}}

\def\tors{_{\textrm{tors}}}      \def\tor{^{\textrm{tor}}}
\def\red{^{\textrm{red}}}         \def\nt{^{\textrm{ssu}}}

\def\sss{^{\textrm{ss}}}          \def\uu{^{\textrm{u}}}
\def\mm{^{\textrm{m}}}
\def\tm{^\times}                  \def\mult{^{\textrm{mult}}}

\def\uss{^{\textrm{ssu}}}         \def\ssu{^{\textrm{ssu}}}
\def\comp{_{\textrm{c}}}
\def\ab{_{\textrm{ab}}}

\def\et{_{\textrm{\'et}}}
\def\nr{_{\textrm{nr}}}

\def\nil{_{\textrm{nil}}}
\def\sol{_{\textrm{sol}}}
\def\End{\textrm{End\,}}

\def\til{\;\widetilde{}\;}


\title[Thompson-like characterization of the solvable radical]
{Thompson-like characterization \\ of the solvable radical}

\author[Guralnick, Kunyavski\u\i , Plotkin, Shalev]{Robert
Guralnick, Boris Kunyavski\u\i , Eugene Plotkin, and Aner Shalev}
\address{Guralnick: Department of Mathematics, University of
Southern California, Los Angeles, CA 90089-2532, USA}
\email{guralnic@usc.edu}
\address{Kunyavski\u\i \ and Plotkin: Department of
Mathematics, Bar-Ilan University, \newline 52900 Ramat Gan,
ISRAEL} \email{kunyav@macs.biu.ac.il, plotkin@macs.biu.ac.il}
\address{Shalev: Einstein Institute of Mathematics,
Edmond J. Safra Campus, Givat
Ram, \newline Hebrew University of Jerusalem, 91904 Jerusalem,
ISRAEL} \email{shalev@math.huji.ac.il}

\begin{abstract}
We prove that the solvable radical of a finite group $G$ coincides
with the set of elements $y$ having the following property: for
any $x\in G$ the subgroup of $G$ generated by $x$ and $y$ is
solvable. We present analogues of this result for finite
dimensional Lie algebras and some classes of infinite groups.
We also consider a similar problem for pairs of elements.
\end{abstract}

\date{\today}
\maketitle

\thispagestyle{empty} 


\section{Introduction} \label{sec:intro}

The first motivating result for this paper is the following famous
theorem of J.~Thompson \cite{Th} (see also \cite{Fl}): a finite
group $G$ is solvable if and only if every 2-generated subgroup of
$G$ is solvable.

Our main aim is to prove the following  extension of Thompson's
theorem.

\begin{theorem} \label{th:finite}
Let $G$ be a finite group, and let $R(G)$ be the solvable radical
of $G$ (namely the maximal solvable normal subgroup of $G$). Then
$R(G)$ coincides with the set of all elements $y\in G$ with the
following property: for any $x\in G$ the subgroup generated by $x$
and $y$ is solvable.
\end{theorem}

Our proof of Theorem 1.1 invokes the Classification of finite
simple groups, and uses the so called ``one and a half
generation'' of almost simple groups, proved by Guralnick and
Kantor \cite{GK} using probabilistic arguments (see Theorem
\ref{3/2} below).
Thus our result may be regarded as yet another demonstration
of the power of probabilistic and counting methods in group
theory (see the survey paper \cite{Sh} for further background).

Theorem 1.1 can be extended to some classes of infinite groups (see
Theorems \ref{th:linear} and \ref{th:PI} below). It has an obvious
Lie-algebraic counterpart (Theorem \ref{th:Lie}). We call a
criterion given in Theorem \ref{th:finite} the Thompson-like
characterization of the solvable radical.

We also strengthen Theorem 1.1, showing that, if $G$ is a finite group,
and $y_1, y_2 \in G$  are two elements both outside $R(G)$, 
then there exists $x \in G$ such that $\langle x, y_i \rangle, i=1,2$ 
are both nonsolvable (see Theorem \ref{th:two} below). 
It is intriguing that this result is best possible, in the sense that 
it cannot be extended for three elements.

Apart from the results, our paper contains various related questions
and conjectures which we plan to consider in the future.


\bigskip

\noindent {\it Acknowledgements}. Guralnick was partially
supported by the National Science Foundation Grant DMS  0140578.
Kunyavski\u\i \ and Plotkin were partially supported by the
Ministry of Absorption (Israel), the Israel Science Foundation
founded by the Israel Academy of Sciences --- Center of Excellence
Program, the Minerva Foundation through the Emmy Noether Research
Institute of Mathematics, and the EU networks HPRN-CT-2002-00287
and INTAS 00-566. Shalev was partially supported by the Israel
Science Foundation founded by the Israel Academy of Sciences.

The authors are very grateful to B. Plotkin and E. Zelmanov for 
inspiring suggestions and discussions.

\section{Lie algebras} \label{sec:Lie}

We start with a Lie-algebraic counterpart of Theorem
\ref{th:finite} which will give us important hints to its proof.

\begin{theorem} \label {th:Lie}
Let $L$ be a finite dimensional Lie algebra defined over a field
$k$ of characteristic zero, and let $R(L)$ be the solvable
radical of $L$ (namely the maximal solvable ideal of $L$). Then $R(L)$
coincides with the set of elements $y\in L$ with the following
property: for any $x\in L$ the subalgebra generated by $x$ and $y$
is solvable.
\end{theorem}

\begin{proof}
If $y\in R(L)$, then for any $x\in L$ the subalgebra generated by
$x$ and $y$ contains a solvable ideal with one-dimensional
quotient and is therefore solvable. We shall give three proofs for
the reverse inclusion.

{\it 1st proof}. Suppose $y\notin R(L)$. We have to prove that
there is $x\in L$ such that the subalgebra generated by $x$ and
$y$ is not solvable. After factoring out $R(L)$, we are reduced to
proving this in the case where $L$ is semisimple. Clearly, it is
enough to consider the case where $L$ is simple. In that case the
result follows from \cite{Io} where it is proved that for any
nonzero element $y$ of a simple Lie algebra $L$ there is $x$ such
that $x$ and $y$ generate $L$.

{\it 2nd proof}. Suppose $y$ has the property stated in the
theorem. We have to prove that $y\in R(L)$. Consider the sequence
of words $v_n(x,y)$ defined inductively by the following rule:
$$
v_1(x,y)=x, v_{n+1}(x,y)=[v_n(x,y),[x,y]],\dots
$$
This sequence can be used for characterization of $R(L)$:
according to \cite[Theorem 3.7]{BBGKP}, $y\in R(L)$ if and only if
for any $x$ there exists $n$ such that $v_n(x,y)=0$.

Let now $x$ be an arbitrary element of $L$. Since the subalgebra
generated by $x$ and $y$ is solvable, it satisfies the identity
$v_n(x,y)\equiv 0$ for some $n$ \cite[Theorem 3.4]{BBGKP}, and we
are done.

{\it 3rd proof}.  Assume that $y \notin R(L)$. We induct on $\dim
L$. Since the set of elements which generate a solvable subalgebra
is closed in the Zariski topology and $L$ is dense in its Zariski
closure (in some algebraic closure of $k$), we may assume that $k$
is algebraically closed. By induction, $L$ is simple (as in the
first proof). If $L$ has rank $1$ (i.e. $L \cong sl(2)$), the
result is straightforward (for any element is in only finitely
many Borel subalgebras and any maximal subalgebra is a Borel). If
the rank of $L$ is larger than $1$,  it is straightforward to see
that there is a maximal parabolic subalgebra $P$ of $L$ such that
$y \in P \setminus R(P)$ (see \cite[Lemma 2.1]{GS} for an easy proof
of this for groups).  Thus, by induction, there exists $x \in P$
with $\langle x, y \rangle$ not solvable.

\end{proof}

\begin{remark} In view of results of Section \ref{sec:lin}, it
seems plausible that Theorem \ref{th:Lie} can be extended to some
classes of infinite dimensional Lie algebras (in particular, to
Lie algebras with polynomial identity).
\end{remark}

Let us consider the case of finite dimensional Lie algebras
of positive characteristic $p$.  The second proof of the theorem
does not extend to this case.   It seems likely that the first
and third proofs may extend. However, the structure of modular Lie
algebras, especially those in characteristics $2$ and $3$, is quite 
a bit more complicated than that of Lie algebras in characteristic 
zero.  We pose the following natural problems:

\begin{quest} \label{modtwo}  Is every simple finite
dimensional Lie algebra generated by two elements?
\end{quest}

If the underlying field is algebraically closed of characteristic
$p > 3$ one can use the available classification to check this
case by case, and a positive answer seems very likely.

The next question concerns one and a half generation.

\begin{quest} \label{modliegen}  Let $L$ be a simple finite 
dimensional Lie algebra.  Is it true that for any $0 \ne y \in L$
there exists $x \in L$ with $L=\langle x,y \rangle$?
\end{quest}


We do point out that, though the problems above are open,
at least solvability of finite dimensional Lie algebras
can be determined by pairs of elements.  
The key to this is the analog of Thompson's theorem
on minimal simple groups (in characteristic zero,
this is quite easy from the basic structure theorems
-- the modular case is quite a bit more difficult, but still
not as difficult as Thompson's result).

We thank A. Premet
for pointing out the following result of Schue \cite{Sc}.
Premet also pointed out that one can also prove this
using \cite{PS}.

\begin{lemma} Let $p > 3$ be a prime and let $F$ be an algebraically
closed field of characteristic $p$.  Let $L$ be a finite dimensional
Lie algebra over $F$ such that every proper subalgebra
is solvable.  Then $L/R(L) \cong sl(2)$, where $R(L)$
is the solvable radical of $L$.
\end{lemma}

There was a mistake in \cite{Sc} but it does  not affect
the statement above.  Also, Schue only concludes that $L/R(L)$
is $3$-dimensional -- but the only simple three dimensional
Lie algebra in characteristic larger than $2$ is $sl(2)$.

This easily yields:

\begin{theorem} \label{modular lie}  Let $p > 3$ be a prime
and let $F$ be an infinite field of characteristic $p$.  Let $L$
be a finite dimensional Lie algebra over $F$ such that
every pair of elements generate a solvable Lie algebra.
Then $L$ is solvable.
\end{theorem}

\begin{proof} Let $E$ be the algebraic closure of $F$
and set $M=EL$. Since the set of pairs $(x,y) \in L \times L$
that generate a solvable Lie algebra is closed (in the
Zariski topology) and since $L \times L$ is dense
in $M \times M$, it follows that $M$ satisfies the same
hypotheses.  So we may assume that $F$ is algebraically
closed.  Induct on $\dim L$.  If $R(L) \ne 0$, we may
pass to the quotient.   So we may assume that $R(L)=0$
and every proper sublalgebra is solvable. By the Lemma,
this implies that $L \cong sl(2)$.  This algebra is clearly
generated by $2$ elements and is not solvable, a contradiction.
\end{proof}

\section{Finite groups} \label{sec:finite}

In this section we prove Theorem \ref{th:finite}. For brevity, let
us introduce the following notion.

\begin{defn} \label{def:rad}
Let $G$ be a group. We say that $y\in G$ is a radical element if
for any $x\in G$ the subgroup generated by $x$ and $y$ is
solvable. Denote by $S(G)$ the set of radical elements of $G$.
\end{defn}

Note that in any group $G$ we have $R(G) \subseteq S(G)$.
Indeed, if $y\in R(G)$, then for any
$x$ the subgroup generated by $x$ and $y$ contains a solvable
normal subgroup with cyclic quotient and is therefore solvable.

With this terminology, Theorem \ref{th:finite} says that if $G$ is
finite, then  $R(G)=S(G)$.

\begin{proof}

It suffices to prove that $S(G) \subseteq R(G)$. One could try to
mimic one of three proofs of Theorem \ref{th:Lie}. Thus,
taking into account the second proof,
Theorem \ref{th:finite} would immediately follow from the
following

\begin{conj} \cite[Conj. 2.12]{BBGKP} \label{conj:BBGKP}
For any finite group $G$ there exists a sequence $u_n(x,y)$ such
that $R(G)$ coincides with the set of elements $y$ having the
following property: for any $x\in G$ there is $n$ such that
$u_n(x,y)=1$.
\end{conj}

However, \cite{BBGKP} contains only partial results towards
Conjecture \ref{conj:BBGKP}.

The reduction argument used in the third proof also encounters
some problems in the group case (say, for groups of Lie type over
small fields).

Fortunately, the ``one and a half generation'' theorem of
Guralnick and Kantor \cite[Corollary of Th.~I on p.~745]{GK}
allows us to imitate the first proof of Theorem \ref{th:Lie}.
We state this theorem as follows.

\begin{theorem} \label{3/2} Let $G$ be a finite almost
simple group with socle $T$. 
Then for every non-identity element $y \in G$
there exists an element $x \in G$ such that  
$\langle x, y \rangle$ contains $T$.
\end{theorem}

First, we use this to prove an auxiliary result.

\begin{lemma} \label{lem:aut}  Let $G$ be a finite group.
Suppose that $N$ is a minimal normal subgroup of $G$ and $G/N$ is
generated by the coset $yN$ with $y$ not in $C_G(N)$.  Then there
exists $x \in N$ with $G=\langle x, y \rangle$.
\end{lemma}

\begin{proof}  We prove this by induction on $|G|$.

If $N$ is abelian, then as $N$ is an irreducible
module for $\langle y \rangle$, $G=\langle x, y \rangle$ for
any nontrivial $x \in N$.

So we may assume that $N=T^n$ where $T$ is a nonabelian simple
group and $y$ transitively permutes the $n$ copies of $T$.  Set
$C:=C_G(N)$. Since $C \cap N=1$ and $N$ is the derived subgroup of
$G$, $C =Z(G)$.  If $C \ne 1$, then by minimality, $G/C = \langle
xC, yC \rangle $ for some $x \in N$. Thus, $\langle x,y \rangle C
$ contains $CN$ and so $\langle x, y \rangle$ contains $N$, the
derived subgroup of $CN$, and so $G=\langle x, y \rangle$.  So by
induction, $C=1$.

We claim that we may assume that $y$ has prime order.
Suppose $y$ has order $pq$ where $p$ is prime and $q > 1$.
Set $w=y^q$ and let $M$ be a normal subgroup of $N$
minimal with respect to being normalized by $w$.  By induction,
$\langle M, w \rangle = \langle x, w \rangle$ for some
$x \in M$. In particular, $\langle x,w \rangle$ contains
a simple direct factor of $N$.  Since $y$ is transitive on the
simple direct factors, it follows that $\langle y, x \rangle$
contains $\langle w, x \rangle$ and so contains all components
of $N$, whence it contains $N$ and $\langle y, x \rangle
=\langle N, y \rangle$ as required.

If $n=1$,
the statement of the lemma
is an immediate corollary of Theorem \ref{3/2}.

If $n > 1$, then $n$ must be prime and $y$ must permute the $n$
copies of $T$.  So we may assume that $y(t_1, \ldots, t_n)y^{-1}
=(t_n,t_1,\ldots, t_{n-1})$.  Choose $t_1, t_n \in T$ such that
$t_n$ is an involution and $T=\langle t_1, t_n \rangle$ (which is
possible by another application of \cite{GK}). Take $x=(t_1,1,
\ldots, 1, t_n)$.  Consider  $J:=\langle x, x^y \rangle$. Let
$\pi$ denote the projection of $N$ onto the first copy of $T$.
Clearly $\pi(J)=T$ (since $\pi(s)=t_1$ and $\pi(s^y)=t_2$).  Also,
$s^2$ is a nontrivial element of $T$ (since $t_n$ is an involution
and $t_1$ is not).   Thus, $[s^2,J]=T \le J$.   Since $y$ permutes
the copies of $T$ transitively, this implies that $N \le \langle
s, y \rangle$, whence $G=\langle s, y \rangle$.

\end{proof}

We are now able to prove that $S(G) \subseteq R(G)$. It is easy to
see that $S(G/R(G)) = S(G)/R(G)$ (more precisely $gR(G) \in
S(G/R(G))$ if and only if $g \in S(G)$). Factoring out $R(G)$, we
may assume that $G$ is semisimple (i.e. $R(G)=1$). We have to
prove that $S(G)=1$.

We use the notion of the generalized Fitting subgroup $F^*(G)$.
See \cite{As} for the basic definitions and results.  The important
facts about $F^*(G)$ for us are that if $R(G)=1$, then
$F^*(G)$ is a direct product of nonabelian simple groups
(the components of $G$) and the centralizer of $F^*(G)$ is trivial.

Since $R(G)=1$, if $1 \ne y \in G$, then there is some component
$L$ of $G$ such that $y$ does not centralize $L$. Let
$\Delta=\{L^g| g \in \langle y \rangle\}$ and set $N$ to be the
subgroup generated by the $L^g \in \Delta$.  Then $N$ is a direct
product of these components. Consider $J:=\langle N, y \rangle$.
Observe that $N$ is the unique minimal normal subgroup of $J$.
Applying
Lemma \ref{lem:aut}
shows that $J=\langle x, y \rangle$ for some $x \in N$. In
particular, $J$ is not solvable and so $S(G)=1$.

The theorem is proved.
\end{proof}

\begin{corol}
Let $G$ be a finite group, let $y \in G$, and let $\left< y^G
\right>$ denote the minimal normal subgroup of $G$ containing $y$.
Then $\left< y^G \right>$ is solvable if and only if the subgroup
$\left<y^{\left< x \right>}\right>$ is solvable for all $x \in G$.
\end{corol}

Since  an element $y\in G$ is radical if and only if the subgroup
$\left<y^{\left< x \right>}\right>$ is solvable for any $x\in G$,
the following question arises:
\begin{quest}\label{rad:el}
For which groups $G$ the normal subgroup $\left< y^G \right>$ is
solvable if and only if the element $y$ is radical?
\end{quest}

\section{Linear groups and PI-groups} \label{sec:lin}

\begin{theorem} \label{th:linear}
Let $K$ be a field. If $G \le \GL (n,K)$, then $R(G)=S(G)$.
\end{theorem}

\begin{proof}
As noted earlier we have $R(G)\subseteq S(G)$, and therefore it
suffices to prove that the subgroup $H$ generated by the set $S(G)$
is solvable (and thus coincides with the radical).

Let $H_1=\left<g_1,\dots ,g_s\right>$ be a finitely generated
subgroup of $H$ where all $g_i$'s are radical elements. Then $H_1$
is approximated by finite linear groups $G_\alpha=H_1/N_\alpha,\
\cap N_\alpha=1$ { in dimension $n$} \cite{Ma}. Each $G_\alpha$ is
finite and is generated by the images of radical elements which are
radical as well, and thus $G_\alpha$ is solvable. Since all
$G_\alpha$'s are linear in dimension $n$, their derived length is
bounded, say, by $k = k(n)$. Thus the group $H_1$ has { derived
length} at most $k$. Each finitely generated subgroup of $H$ lies in
some $H_1$. Thus $H$ is locally solvable. Since $H$ is linear, it is
solvable \cite{Za}.
\end{proof}


For the case of PI-groups we use some facts from \cite{Pi},
\cite{Plo1}, \cite{Plo}.

\begin{defn} \label{def:pi} A group $G$ is called a PI-group
(PI-representable in terms of \cite{Pi}) if $G$ is a subgroup of
the group of invertible elements of an associative PI-algebra over
a field.
\end{defn}

Linear groups are a particular case of PI-groups. It is known that
every PI-group $G$ has a unique maximal locally solvable normal
subgroup $R(G)$ called the locally solvable radical of $G$ and that
the locally solvable radical of a finitely generated PI-group is
solvable \cite{Pi}. (For arbitrary groups the locally solvable
radical may not exist, and for arbitrary PI- groups the locally
solvable radical is not necessarily solvable).

PI-groups have the following invariant series: $1\triangleleft
H_0\triangleleft H\triangleleft G$ where $H_0$ is a locally
nilpotent normal subgroup, $H/H_0$ is nilpotent and $G/H$ is a
linear group over a cartesian sum of fields \cite{Plo}.

We want to show that in a PI-group $G$ the locally solvable
radical coincides with $S(G)$.


Let us introduce a useful notion of oversolvable group.

\begin{defn} (cf. \cite{Plo1}). A group $G$ is called oversolvable
if it has an ascending normal series with locally nilpotent factors.
\end{defn}

An arbitrary group $G$ has the oversolvable radical $\widetilde
{HP}(G)=\widetilde \eta(G)$ (that is the unique maximal normal
oversolvable subgroup). The quotient
group $G/\widetilde {\eta}(G)$ is semisimple with respect to the
property of being locally nilpotent, i.e. $\eta(G/\widetilde
{\eta}(G))=1$ where $\eta(G)$ is the locally nilpotent radical of
$G$ (see \cite{Plo1} for the above facts). If $G$ is finite, noetherian,
or linear, $\widetilde {\eta}(G)$ coincides with the solvable
radical $R(G)$ \cite{Su}.

\begin{theorem} \label{th:PI}
If $G$ is a PI-group, then $R(G)=\widetilde{\eta}(G)= S(G)$.
\end{theorem}

\begin{proof}
We consider three cases. 1. $G \le \GL_n(P)$ where $P$ is a field.
2. $G \le \GL_n(K)$ where $K$ is a cartesian sum of fields. 3.
General case.

Case 1. If $G \le \GL_n(P)$, where $P$ is a field, then $R(G)=S(G)$
by Theorem \ref{th:linear}.

Case 2. Suppose $G \le \GL_n(K)$, where $K=\bigoplus_s P_s$ is a
cartesian sum of fields.  Consider the set of congruence subgroups
$U_s$ of $\GL_n(K)$ such that $\GL_n(K)/U_s\cong\GL_n(P_s)$. Since
$\bigcap_s U_s=1$, the group $G$ lies in the cartesian product
$\prod_s\GL_n(P_s).$ We have
$G\subset\GL_n(K)\subset\prod_s\GL_n(P_s)\to \GL_n(P_s) $.
Let $H$ be the subgroup in $G$ generated by the set $S(G)$. It is enough
to show that $H$ is solvable. Set $U_s'=H\cap U_s$. Then $\bigcap_s U_s'=1$.
Each $H/U_s'$ can be viewed as a subgroup in $\GL_n(P_s)$ and is therefore
generated by the radical elements. Thus they are all solvable of bounded
derived length. Therefore $H$ is solvable and $R(G)=S(G)$.
Similar arguments give $R(G)=\widetilde {\eta}(G)$.

Case 3. Let us first show that $R(G)\subseteq S(G)$. Let $g\in
R(G)$, $h\in G$. Consider the subgroup $G_0=\left<g, h\right>$. We
have $g\in R(G)\cap G_0$ and, consequently, $g\in R(G_0)$. By
\cite{Pi}, the locally solvable radical of a 2-generated group is
solvable. So the group $G_0$ is solvable as a cyclic extension of
a solvable group.

Now we want to show that $R(G)=\widetilde{\eta}(G)$. Let us first
prove that $\widetilde{\eta}(G)\subseteq R(G).$ We have to prove
that the group $\widetilde{\eta}(G)$ is locally solvable. We take a
finitely generated subgroup $G_0$ in $\widetilde{\eta}(G)$ and
show that $G_0$ is solvable. Consider the locally solvable radical
$R(G_0)$. Since $G_0$ is finitely generated, the radical $R(G_0)$
is solvable \cite{Pi}. So it is enough to prove that the group
$G_0/R(G_0)$ is solvable. We use the following result about the
structure of PI-groups \cite{Plo}: in every PI-group $G$ the
quotient group $\widetilde{\eta}(G)/ \eta(G)$ is solvable. We have
$\eta(G)\subseteq R(G)$ as a locally nilpotent subgroup.  Apply
this to $G_0$. Since $G_0\subseteq\widetilde{\eta}(G)$, we have
$\widetilde{\eta}(G_0)=G_0$. Since $\widetilde{\eta}(G_0)/
\eta(G_0)$ is solvable and $\eta(G_0)\subseteq R(G_0)$, the group
$\widetilde{\eta}(G)/R(G_0)$ is solvable. Thus $G_0/R(G_0)$ is
solvable, and hence so is $G_0$. The inclusion
$\widetilde{\eta}(G)\subseteq R(G)$ is proved.

Let us prove the opposite inclusion
$R(G)\subseteq\widetilde{\eta}(G)$. Recall that in every PI-group
$G$ there is a normal subgroup $H$ which is an extension of a
locally nilpotent group by a nilpotent group and such that $G/H$
lies in $\GL_n(K)$ where $K$ is a cartesian sum of fields. Then
$H$ is oversolvable. Therefore $H\subseteq \widetilde {\eta}(G)$,
and thus $H\subseteq R(G)$. Consider the group $G/H$ and its
subgroup $R(G)/H$. This is a locally solvable normal subgroup in
$G/H$ and thus lies in $R(G/H)$. The group $G/H$ is linear and
hence $R(G/H)= \widetilde{\eta}(G/H)$. Then $R(G)/H\subseteq
\widetilde{\eta}(G)/H=\widetilde{\eta}(G/H) $. Thus
$R(G)\subseteq\widetilde{\eta}(G)$ and $R(G)=\widetilde{\eta}(G)$.

We are now able to prove that $S(G)\subseteq R(G)$. Let $g\in
S(G)$. Denote by $\bar g\in G/H$ its image under the natural
projection. Then $\bar g\in R(G/H)$ and thus $\bar
g\in\widetilde{\eta}(G)/H. $ Then $g\in\widetilde{\eta}(G)=R(G)$.
\end{proof}

The above theorem has an obvious consequence which can be viewed
as a natural generalization of Thompson's theorem:

\begin{corol}
A PI-group $G$ is locally solvable if and only if every
two-generated subgroup of $G$ is solvable.
\end{corol}

\begin{remark}
In linear groups the locally solvable radical is solvable. From the
above theorem it follows that in PI-groups the locally solvable
radical is solvable modulo the locally nilpotent radical. Indeed,
$\widetilde{\eta}(G)/ \eta(G)$ is solvable \cite{Plo}, and
$\widetilde{\eta}(G)=R(G)$ by Theorem \ref{th:PI}.
\end{remark}

\begin{corol} If $G$ is a PI-group then the normal subgroup
$\left< y^G \right>$ is locally solvable if and only if the
element $y$ is radical. In particular, if $G$ is a finitely
generated PI-group or a linear group, then  $\left< y^G \right>$
is solvable if  and only if the element $y$ is radical.
\end{corol}

\section{Residually finite groups and Burnside-type problems}

We extend some of the results to residually finite groups (in
particular to profinite groups). There are some open questions.  A
trivial consequence of
Theorem \ref{th:finite}
is the following:

\begin{corol}  Let $G$ be a residually finite group.  Then
$\langle y^G \rangle $ is residually finite solvable if and only
if $\langle x, y \rangle$ is residually finite solvable for all $x
\in G$.
\end{corol}

Regarding Question \ref{rad:el} (Section 3), if we do not assume
that the group $G$ is finitely generated, it is easy to construct
counterexamples.

For instance, let $G_i$ be a finite solvable group of derived length
$i$ generated by two elements $x_i, y_i$. Let $P$ be the
direct product of the $G_i$ and $S < P$ the direct sum.  Let
$y=(y_1, y_2, \ldots ) \in P\setminus S$ and set $G=\langle S, y
\rangle$. Clearly, $P$ (and so $G$) is residually finite. Suppose
that $x \in G$.  Then $x=sy^j$ for some $s \in S$ and integer $j$.
Setting $H= \langle x, y \rangle$, we see that $[H,H]$ is
contained in a finite direct product of the $G_i$ and so is
solvable.  Hence so is $H$, and thus $y$ is a radical element of
$G$. On the other hand, the normal closure of $y$ in $G$ is
not solvable (since it contains $[G,G]$ and $G$ is not solvable).

Moreover, for finitely generated residually finite groups (or
profinite groups) the statement of Question \ref{rad:el} is false
as well.

\begin{prop}
There exists a finitely generated residually finite group $G$ such
that for some radical element $y\in G$ the group $\langle y^G
\rangle $ is not solvable.
\end{prop}

\begin{proof}
Let $G$ be a three generated residually finite group that is not
nilpotent and in which every two generated subgroup is a nilpotent
group. Such groups exist due to Golod--Shafarevich \cite{Go}.

We show that the group $G$ provides a counterexample to the
statement of Question \ref{rad:el}.

Indeed, every element of $G$ is, obviously, Engel and radical.
Take any element $y\in G$ which does not belong to the locally
nilpotent radical of $G$. We claim that the normal
subgroup $\left< y^G \right>$ is not solvable. For a solvable
normal subgroup consisting of Engel elements should be locally
nilpotent, and thus should be contained in the locally nilpotent radical,
contradicting the choice of $y$. 
We conclude that $y$ is a radical element of $G$, but its normal
closure $\left< y^G \right>$ is not solvable.
\end{proof}

Note that we have a stronger necessary condition for $y^G$ to
generate a solvable group.  If $\langle y^G \rangle$ is solvable
of derived length $d$, then $\langle x, y \rangle$ is solvable of
derived length at most $d+1$ for every $x \in G$.

\begin{quest} Let $G$ be a residually finite group and suppose
$y \in G$ satisfies
$\langle x, y \rangle$ is solvable of derived length at most $d$
(for some positive integer $d$ independent of $x$). Does it follow
that $\langle y^G \rangle$ is solvable of derived length at most 
$f(d)$ for some function $f$?
\end{quest}

Note that this reduces to solving the problem for finite groups.
An affirmative answer would give a characterization of the set of
elements in a residually finite group whose normal closure is
solvable.

The Golod example shows that every $2$-generated subgroup being
solvable does not imply that $G$ is solvable for $G$ a finitely
generated residually finite group. So we ask:

\begin{quest} Let $G$ be a residually finite group and suppose
$\langle x, y \rangle$
is solvable of derived length at most $d$ for every pair of elements 
$x,y \in G$. Does it follow that $G$ is solvable of derived length 
at most $f(d)$ for some function $f$?
\end{quest}

This reduces to considering finite solvable groups.
It is also natural to consider the following variation, where we bound
the number of generators of $G$.

\begin{quest} Let $G$ be a residually finite group generated by
$c$ elements, and suppose $\langle x, y \rangle$
is solvable of derived length at most $d$ for every pair of elements 
$x,y \in G$. Does it follow that $G$ is solvable of derived length at most
$f(c,d)$ for some function $f$?
\end{quest}

Again this reduces to considering finite solvable groups.
In fact, some of the questions above may be even posed in
greater generality, namely for arbitrary groups.
The obvious analogous problems for Lie algebras are also of 
interest. In particular:

\begin{quest} Let $L$ be a Lie algebra in which every two elements
generate a solvable subalgebra of derived length at most $d$.
Does it follow that $L$ is solvable? or locally solvable?
\end{quest}

The problems posed in this section may be regarded as
Burnside-type problems; it would be interesting to find
out whether Burnside-type techniques can help tackling them.

\section{Pairs of elements}

We can extend some of these results to pairs of elements not in
the solvable radical.

We first state an even stronger version for Lie algebras.

\begin{theorem}  Let $L$ be a finite dimensional Lie algebra
over a field $k$ of characteristic zero.  Let $X$ be a finite
set of
elements of $L \setminus{R(L)}$.  Then there is some $y \in L$
such that $\langle x, y \rangle$ is not solvable for all $x \in
X$.
\end{theorem}

\begin{proof}  Fix $x \in X$ and set $L(x):=\{y \in L|\langle x, y \rangle
\ \mathrm{is \ solvable} \}$.  Note that $L(x)$ is a closed subvariety
of $L$ (in the Zariski topology -- for the solvable condition is given
by certain words in $x$ and $y$ being trivial).
By Theorem \ref{th:Lie}, $L(x)$
is a proper subvariety.  Affine space over an
infinite field cannot be written
as a finite union of proper closed subvarieties and
$\cup_{x \in X} L(x) \ne L$,
whence the result.
\end{proof}

We have a similar result for connected algebraic groups (and for
connected Lie groups).

\begin{theorem}
Let $G$ be a connected algebraic group over an
algebraically closed field $k$.  Let $X$ be any finite set of
elements outside $R(G)$.  Then there exists $y \in G(k)$ such
that $\langle x, y \rangle$ is not solvable for all $x \in X$.
\end{theorem}

The proof is identical (for a fixed $x \in G(k)$, the set of $y
\in G(k)$ with $\langle x, y \rangle$ solvable is a closed
subvariety; if $x$ is not in $R(G)$, it is proper by earlier
results and so the union is also proper).

\medskip

We have a weaker (but harder) result for finite groups.  Clearly,
we cannot take any finite subset as above (for example,
if $X$ is the set of all nontrivial elements of $G$).
See \cite{BGK} for references on this general problem.

Here is an example to show that even for subsets of
size $3$, there can be a problem.  Consider $G=A_5$.  Let
$x_1=(23)(45)$, $x_2=(13)(45)$ and $x_3=(12)(45)$.
Note that no two of the $x_i$ normalize a common
Sylow $5$-subgroup (since the product of any two distinct
$x_i$ has order $3$ and the normalizer of a Sylow $5$-subgroup
is dihedral of order $10$).  Since there are six Sylow $5$-subgroups,
it follows that if $y$ has order $5$, $\langle x_i, y \rangle$
is dihedral of order $10$
for exactly one $x_i$.  If $y$ has order $3$, then
either $y$ fixes $i$ for some $i \le 3$
and so $\langle x_i, y \rangle$ is contained
in the stabilizer of $i$ (i.e.  $A_4$)
or $y$ fixes $4$ and $5$ and so $\langle x_i, y \rangle$
is contained in an $S_3$ for each $i$.  If $y$ has order $2$,
then $\langle x_i,y \rangle$ is a dihedral group for each $i$.
So we have shown that for any $y \in G$, for at least one
$i$,  $\langle x_i,y \rangle$ is solvable.

In order to prove the analogous result for finite groups (with
$|X|=2$), we need Theorem 1.4 from \cite{BGK}. This extends the
result of \cite{GK} (see Theorem \ref{3/2} above).
Both results are proved using probabilistic methods.

\begin{theorem} \label{BGK1.4} {\rm{\cite[Theorem 1.4]{BGK}}}
Let $G$ be a finite almost simple group with socle $S$.
If $x,y$ are nontrivial elements of $G$, then there exists
$s \in S$ such that $S$ is contained in $\langle x, s \rangle$
and $\langle y, s \rangle$.
\end{theorem}

The key point that we will use is that the subgroups
$\langle x, s \rangle$ and $\langle y, s \rangle$
are not solvable.

We can now prove the main result of this section.

\begin{theorem}\label{th:two}  Let $G$ be a finite group.  Suppose that
$x$ and $y$ are not in $R(G)$.  Then there exists $s \in G$
such that $\langle x, s \rangle$ and $\langle y, s \rangle$
are not solvable.
\end{theorem}

\begin{proof}  Suppose that $G$ is a counterexample
of minimal order. Thus, $R(G)=1$ (or we could
pass to $G/R(G)$). We can
replace $x$ and $y$ by powers and so assume that they
each have prime order (this is not essential to the
proof given below -- it does make it a bit easier
to see the possibilities for the action of $x$ or $y$).

Then $G$ has a normal subgroup $N$ (the generalized Fitting
subgroup) that is the direct product of simple groups $L_i, 1 \le
i \le n$, and $C_G(N)=1$. In particular, this implies that
$\langle x^N \rangle$ is not solvable (for since $x$ does not
centralize $N$, $\langle x^N \rangle$ must intersect $N$ in a
nontrivial normal subgroup and any normal subgroup of $N$ is a
direct product of simple groups and in particular is not
solvable). So already the hypotheses are satisfied in 
$\langle x, y, N \rangle$ and so the minimality hypothesis implies that
$G=\langle x, y, N \rangle$.

The direct factors of $N$ are the components of $G$.
First suppose that there is some component $L$ of $G$
with neither $x$ nor $y$ in $C_G(L)$.  Let $M$ be the
normal closure of $L$ in $G$.  Then we may assume
that $M=N$ (otherwise, $C_N(M) \ne 1$ and we can
pass to the smaller group $G/C_N(M)$ and arguing as
above, we see that the normal closures of $x$ and $y$
are still not solvable).

Write $N=L_1 \times \ldots \times L_m$ with say $L=L_1$. If $m=1$,
then $G$ is almost simple and we apply Theorem \ref{BGK1.4} to
obtain the conclusion.

So assume that $m > 1$ and that $L_1^x = L_2$.
If $L_1^y \ne L_1$ or $L_2$, choose $u,v \in L_1$
that generate $L_1$ and consider the element
$s=uv^xv^y$.  Note that $\langle s, s^{x^{-1}} \rangle$
projects onto $L$ and so $\langle s, x \rangle $ is not
solvable and similarly for $\langle s, y \rangle$.

If $L_1^y=L_2$, choose $u,v \in L$ such that $\langle u, v \rangle
= \langle u, v^{xy^{-1}} \rangle = L$ (this can be done by
\cite[Theorem 1.2 or 1.4]{BGK}).  Then $s=uv^{x}$ satisfies the
conclusion.

Finally suppose that $L_1^y=L_1$.  So $y$ induces a nontrivial
automorphism of $L_1$. Then choose $u,v \in L$ (by Theorem
\ref{BGK1.4}) such that $L = \langle u, v \rangle = \langle
u^{\langle y \rangle}, i = 0, 1, \ldots \rangle$ and again set
$s=uv^{x}$.

The remaining case to consider is that for each component $K$ of
$G$, either $x$ or $y$ centralizes $K$ (but not both by hypothesis).

So choose components $L$ and $K$ with $x$ not centralizing $L$ and
$y$ not centralizing $K$.   Arguing as above, we can choose $a\in
\langle L, L^x \rangle$ with $H_x:=\left<x,a\right>$ not solvable
and $b \in \langle K, K^{y} \rangle$ with $H_y:=\langle y,b
\rangle$ not solvable. Set $s=ab$. Note that $ay=ya$ and $xb=bx$.
Let $D:=\langle s, x \rangle$. Then $D$ centralizes $b$ and so
$\langle D, b \rangle =D \langle b \rangle$ contains $H_x$ and is
not solvable.  Thus, $D$ itself is not solvable.  Similarly,
$\langle s, y \rangle$ is not solvable and the theorem is proved.
\end{proof}

\begin{corollary}  Let $G$ be a linear group over a field
$k$.  If $x$ and $y$ are not in the
solvable radical of $G$, then there exists $s \in G$ with
$\langle x, s \rangle$ and $\langle y, s \rangle$ not solvable.
\end{corollary}

\begin{proof} We may assume that $G$ is finitely generated. By our main
result on linear groups,  there is a finite homomorphic image $H$
of $G$ in which the images of $\langle
x^G \rangle$ and $\langle y^G \rangle$ are not solvable.
So we may choose $\bar{s} \in H$ such that $\langle \bar{x}, \bar{s} \rangle$
and $\langle \bar{y}, \bar{s} \rangle$ are not solvable.  Now take $s$ to
be any preimage of $\bar{s}$.
\end{proof}


\section{Concluding remarks}

Let us observe that certain important classes of groups and Lie
algebras can be explicitly characterized in terms of two-variable
identities: one can mention here classical results for finite
dimensional nilpotent Lie algebras (Engel) and finite nilpotent
groups (Zorn \cite{Zo}) and their recently obtained counterparts for
finite dimensional solvable Lie algebras \cite{GKNP} and finite (or
linear) solvable groups \cite{BGGKPP} (see also \cite{BWW}).
Moreover, Engel identities were used by Baer to characterize
explicitly the nilpotent radical of an arbitrary finite (and, more
generally, noetherian) group \cite{Ba}. Baer's theorem was extended
to the locally nilpotent radical of linear groups and PI-groups
\cite{Pla}, \cite{Plo} and to the nilpotent and the solvable radical
of finite dimensional Lie algebras \cite{BBGKP}. These results give
a certain hope for characterization of the solvable radical $R(G)$
of a finite group $G$ in similar, Engel-like terms.
However, the corresponding Conjecture 3.2 (see also \cite{BBGKP}) 
is still far from being proved, and therefore less explicit
descriptions of the solvable radical, such as the Thompson-like
characterization of Theorem \ref{th:finite}, are very useful.





\end{document}